# Two Precision-controlled Numerical Algorithms for the CDF of Doubly Non-central Beta Distribution Based on the Segmentation of the Infinite Double Series Matrix

**Han Li [1], Fangfang Ma [1,*], Junjie Wang [2,*], Yinhua Tian [1], Baoli Dai [1] and Tianyan Dong [1]**

[1] Shandong University of Science and Technology; lihan_60@sdust.edu.cn, skd993125@sdust.edu.cn, skdxxtyh@163.com, 202123040102@sdust.edu.cn, dty02@sdust.edu.cn
[2] Taishan Institute of Science and Technology; wangjunjie@ta.shandong.cn
* Correspondence: skd993125@sdust.edu.cn, wangjunjie@ta.shandong.cn; Tel.: +86-1521-538-7676, +86-1361-538-7191

**Abstract:** The cumulative distribution function (CDF) of the doubly non-central beta distribution can be expressed as an infinite double series. By truncating the sum of this series, one can obtain an approximate value of the CDF. Although numerous methods exist for calculating the non-central beta distribution, which allow for the control of the truncation range and estimation of the computational error, no such methods have been developed for the doubly non-central beta distribution. In this paper, we propose two new numerical computation methods based on the segmentation of the infinite double series, termed DIV1 and DIV2. Both methods enable automated calculations once the error control parameters are set; there is no need to predetermine the truncation range, and their computational times are comparable. Following detailed derivations, we have established the upper bounds of the errors for both methods, thus ensuring the determinability of the precision.

## 1. Introduction

If $X_1$ and $X_2$ are independent non-central $\chi^2$ random variables with $n_1$ and $n_2$ degrees of freedom and non-centrality parameters $\lambda_1$ and $\lambda_2$, respectively, then the ratio $B = X_1/(X_1 + X_2)$ is referred to as the doubly non-central beta distribution. This distribution has shape parameters $n_1/2$ and $n_2/2$, and non-centrality parameters $\lambda_1$ and $\lambda_2$. Let $F = (X_1/n_1)/(X_2/n_2)$. The random variable $F$ follows a doubly non-central $F$ distribution with $n_1$ and $n_2$ degrees of freedom, and non-centrality parameters $\lambda_1$ and $\lambda_2$. It can be easily shown that $B = n_1F/(n_1F + n_2)$. Therefore, the doubly non-central beta and $F$ distributions can be related and computed from one another [1].

The doubly non-central beta and F distributions are widely utilized in practice. Li et al. [2] confirmed that when a signal is present, the detection statistic, namely the power spectrum sub-band energy ratio, follows the doubly non-central beta distribution. Feng [3,4] employed the non-central F distribution to develop an integrity monitoring system for carrier phase ambiguities. To assess signal transmission quality, Jeske and Sampath [5] proposed an improved signal related to the interference plus noise ratio, which adheres to the singly non-central F distribution, a special case of the doubly non-central F distribution. Preisig and Johnson [6] introduced a doubly non-central F distribution statistic for underwater sonar signal detection. In recent years, the multivariate coefficient of variation (MCV) has garnered increasing attention from researchers as an effective tool for process monitoring in statistical process control [7–10]. Notably, MCV follows the singly non-central F distribution [11]. Ayyoub et al. [12] utilized the doubly non-central F distribution for monitoring the multivariate coefficient of variation. Importantly, in these applications, it is essential to calculate the CDF for both the doubly non-central beta and F distributions.

The CDF of the doubly non-central beta and F distributions can be expressed using an infinite double series [13], which serves as the foundation for various calculation methods. Tiku [14] utilized three moments to approximate the CDF of the doubly non-central F distribution. While this method allows for rapid computation, it is generally only accurate to the third decimal place. When both $\lambda_1$ and $\lambda_2$ are large, or when one is large and the other is small, the results are expected to exhibit significant error. Tiku and Yip [15] considered four moments to approximate the CDF of the doubly non-central F distribution, but this form is accurate to only the fourth decimal place. Both Tiku [16] and Chattamvelli [17] used Laguerre polynomials to express doubly non-central beta and F distributions, but the application scope of this method is restricted. These methods involve two problems: the calculation range cannot be automatically defined, and the calculation precision cannot be determined. There are two main reasons for these problems. First, the CDF of the abovementioned representations can only be computed numerically, owing to which the number of items in the series is limited. In fact, the number of items is selected according to the results of numerical experiments and experience of researchers. If the number of items is excessively small, the calculation result will be inaccurate; otherwise, an excessive calculation time is required. Second, the calculation errors associated with the abovementioned methods are difficult to determine, and thus, the calculation precision cannot be ensured. In particular, when the number of calculated items is large, it is extremely difficult to choose a reasonable number within a pre-determined error.

Although there is no method for calculating doubly non-central beta distribution that can control the truncating range and estimate the calculation error, but there are still many methods for calculating non-central beta or F distribution. Because of the complexity of series calculation, it is difficult to find a method to determine the calculation error directly. But the error can be controlled indirectly through the upper error bound. Norton [18] derived a error bound to control the calculating error indirectly. Lenth [19] followed Norton's method [18] and derived another error bound. Referring to the method of controlling the calculating error indirectly, this paper proposes two methods for automatically calculating the CDF of doubly non-central beta and F distributions according to the preset upper error bound.

## 2. Theories and existing methods

### 2.1. Doubly non-central beta distribution

According to the definition of the G distribution and doubly non-central beta distribution, presented in Chapter 30 of the book of Johnson et al. [1], the CDF of the doubly non-central beta distribution is

$$B_{n_1,\,n_2}^{\lambda_1,\,\lambda_2}(x) = e^{-\frac{1}{2}(\lambda_1+\lambda_2)} \sum_{j=0}^{\infty} \sum_{l=0}^{\infty} \frac{\left(\frac{\lambda_1}{2}\right)^j \left(\frac{\lambda_2}{2}\right)^l}{j!l!} I_x(j+\frac{n_1}{2}, l+\frac{n_2}{2}), \quad 0 < x < 1. \tag{1}$$

where $I_x(a,b)$ is the incomplete beta function:

$$I_x(a,b) = \frac{\int_0^x t^{a-1}(1-t)^{b-1}dt}{B(a,b)} \leq 1,$$

where $B(a,b) = \Gamma(a)\Gamma(b)/\Gamma(a+b)$. $I_x(a,b)$ has two important properties that will be used in later derivations, as follows:

**Theorem 1.** *When* $0 \leq x \leq 1$, $a, b \in N_+$, *then* $I_x(a-1,b) > I_x(a,b)$ .

**Proof of Theorem 1.** According to the equation 26.5.16 in Abramowitz [20], the following expression can be obtained:

$$I_x(a-1,b) = \frac{\Gamma(a+b-1)}{\Gamma(a)\Gamma(b)} x^{a-1}(1-x)^b + I_x(a,b) \quad . \tag{2}$$

Obviously, $x^{a-1}(1-x)^b > 0$, therefore $I_x(a-1,b) > I_x(a,b)$ . □

**Theorem 2.** *When* $0 \leq x \leq 1$, $a, b \in N_+$, *then* $I_x(a,b-1) < I_x(a,b)$ .

**Proof of Theorem 2.** The equation 26.5.10 in [20] is as follows:

$$I_x(a,b) = xI_x(a-1,b) + (1-x)I_x(a,b-1) \quad . \tag{3}$$

By substituting Equation (2) into Equation (3), the following equation can be obtained:

$$I_x(a,b) = \frac{\Gamma(a+b-1)}{\Gamma(a)\Gamma(b)} x^a(1-x)^b + xI_x(a,b) + (1-x)I_x(a,b-1) \quad .$$

This equation can be simplified as

$$I_x(a,b) = \frac{\Gamma(a+b-1)}{\Gamma(a)\Gamma(b)} x^a(1-x)^{b-1} + I_x(a,b-1) \quad .$$

Because $x^a(1-x)^{b-1} > 0$, therefore $I_x(a,b-1) < I_x(a,b)$. □

To simplify the algebraic representations, let $a = n_1/2$, $b = n_2/2$, $\delta_1 = \lambda_1/2$, and $\delta_2 = \lambda_2/2$. Equation (1) can be rewritten as

$$B_{a,\,b}^{\delta_1,\,\delta_2}(x) = e^{-(\delta_1+\delta_2)} \sum_{j=0}^{\infty} \sum_{l=0}^{\infty} \frac{\delta_1^j \delta_2^l}{j!l!} I_x(j+a, l+b) \quad , \quad 0 < x < 1 \ . \tag{4}$$

In the process of calculating the CDF of the doubly non-central beta distribution, the calculation of the incomplete beta function is extremely complex and requires extensive computing time. When at least one of the degrees of freedom $n_1$ or $n_2$ is an even value, Singh and Relyea [21] obtained the exact expression for $I_x(a,b)$ as follows:

$$I_x(a,b) = \begin{cases} 1-(1-x)^b\left(1+\sum\limits_{j=1}^{a-1}\left(\prod\limits_{i=1}^{j}\dfrac{b+j-1}{i}\right)x^j\right) & ,\ n_1 \textit{ is an even integer} \\ x^a\left(1+\sum\limits_{j=1}^{b-1}\left(\prod\limits_{i=1}^{j}\dfrac{a+i-1}{i}\right)(1-x)^j\right) & ,\ n_2 \textit{ is an even integer} \end{cases} \tag{5}$$

When $n_1$ and $n_2$ are odd integers, Singh and Relyea [21] presented the calculation of $I_x(a,b)$ as follows:

$$\begin{aligned} I_x(a,b) &= \frac{1}{2} + \frac{2}{\pi}(x(1-x))^{1/2}\, D1(D2-D3) - \frac{1}{\pi}\sin^{-1}(1-2x) \ , \textit{where} \\ D1 &= \prod_{j=1}^{a-.5}\left(\frac{jx}{j-.5}\right) \ , \quad D2 = \left(\sum_{k=1}^{b-.5}\left(\prod_{j=1}^{k-1}\frac{a+j-.5}{j+.5}\right)\right)(1-x)^{k-1}, \ \textit{and} \\ D3 &= \sum_{k=1}^{a-.5}\left(\prod_{j=1}^{k-1}\frac{j}{j+.5}\right)x^{k-1} \quad . \end{aligned} \tag{6}$$

As indicated in Equations (5) and (6), the calculation of the incomplete beta function is a time-consuming task in the numerical computation.

Equation (4) is represented by an infinite double series that has no analytical expression. Its results can only be obtained through a numerical computation. The CDF of the doubly non-central F distribution can be obtained from that of the doubly non-central beta distribution as follows:

$$F_{a,\,b}^{\delta_1,\,\delta_2}(f)=\Pr\left(\frac{X_1 n_2}{X_2 n_1}<f\right)$$
$$=\Pr\left(\frac{X_1}{X_1+X_2}<\frac{n_1 f}{n_1 f+n_2}\right)=B_{a,\,b}^{\delta_1,\,\delta_2}\left(\frac{n_1 f}{n_1 f+n_2}\right) \tag{7}$$

The series item of $B_{a,\,b}^{\delta_1,\,\delta_2}(x)$ be expressed as

$$L_{j,\,l}(x;\delta_1,\delta_2,a,b)=e^{-(\delta_1+\delta_2)}\frac{\delta_1^j\delta_2^l}{j!l!}I_x(j+a,l+b) \quad .$$

To simplify the algebraic representations, $L_{j,\,l}(x;\delta_1,\delta_2,a,b)$ can be abbreviated as $L_{j,\,l}$. $B_{a,\,b}^{\delta_1,\,\delta_2}(x)$ can be rewritten as

$$B_{a,\,b}^{\delta_1,\,\delta_2}(x)=\sum_{j=0}^{\infty}\sum_{l=0}^{\infty}L_{j,\,l} \quad .$$

If $j$ and $l$ are considered as the row and column indices, respectively. Then, matrix $\mathbf{M}$ can be defined as follows:

$$\mathbf{M}=\begin{bmatrix} L_{0,\,0} & L_{0,1} & L_{0,\,2} & \cdots \\ L_{1,\,0} & L_{1,1} & L_{1,\,2} & \cdots \\ L_{2,\,0} & L_{2,1} & L_{2,\,2} & \cdots \\ \vdots & \vdots & \vdots & \vdots \end{bmatrix}.$$

In fact, $B_{a,\,b}^{\delta_1,\,\delta_2}(x)$ is the sum of all the elements in the $\mathbf{M}$-matrix. When $a=20$, $b=492$, $x=0.1$, $\delta_1=30.72$, and $\delta_2=20.48$, the item values of $\mathbf{M}$ are as shown in Figure 1.

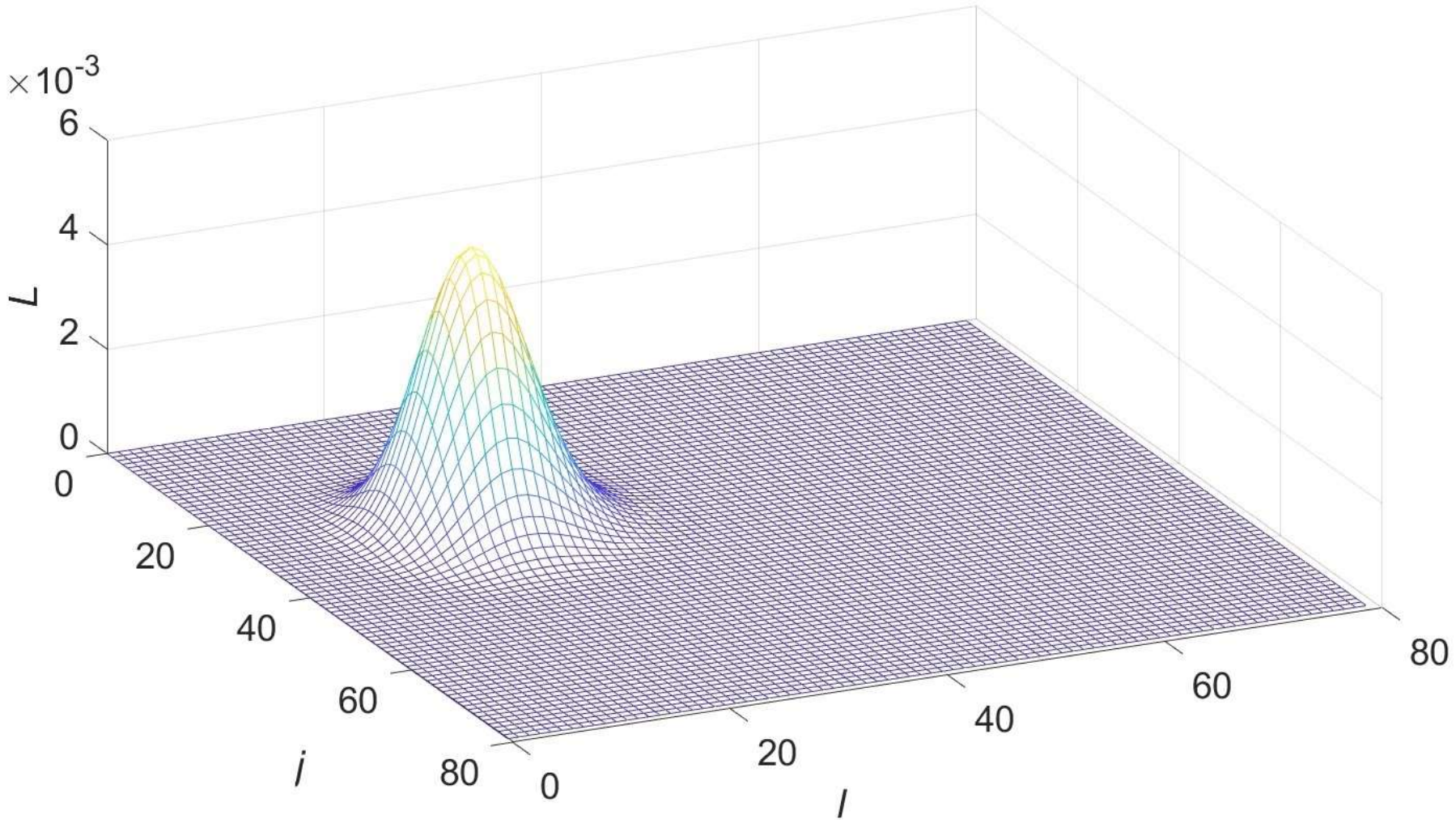

**Figure 1.** Item values of the $\mathbf{M}$-matrix with a=20, b=492, x=0.1, $\delta_1=30.72$, and $\delta_2=20.48$.

In Figure 1, the upper left part of the $\mathbf{M}$-matrix has a notable convex region. The value in the region is large. This indicates that the contribution of the convex region to $B_{a,\,b}^{\delta_1,\,\delta_2}(x)$ is significant. Therefore, as long as the values of the convex region are included in the cumulative calculation process, the calculating precision of the approximate value for CDF can be ensured.

### 2.2. Existing algorithms

The CDF of the doubly non-central beta distribution can be accurately calculated using two class approaches based on direct calculation or Laguerre polynomial calculation, respectively.

The method of direct calculation can be directly calculated according to Equation (4). However, the truncating range should be set in advance before the calculation. Therefor, the number of rows and columns in the upper left part of the **M**-matrix should be determined. As long as the truncating range is sufficiently large, a high calculation precision can be ensured. Bulgren [13] used a direct calculation to obtain the CDF of several doubly non-central F distributions, with precision up to the fourth decimal place.

Another method to calculate the CDF of a doubly non-central beta distribution is through Laguerre polynomials. When $b$ having an integer value, Chattamvelli [17] derived the CDF by Laguerre polynomials as follows:

$$B_{a,\,b}^{\delta_1,\,\delta_2}(x)=e^{-\delta_1(1-x)}x^a\sum_{j=0}^{\infty}\frac{e^{-\delta_2}\delta_2^l}{j!}\sum_{k=0}^{b+j-1}(1-x)^k L_k^{(a-1)}(-\delta_1 x) \ ,$$

where $L_k^{(a-1)}(-\delta_1 x)$ represents Laguerre polynomials, defined as

$$L_k^{(a-1)}(-\delta_1 x)=\sum_{i=0}^{j}\binom{j+a-1}{j-i}\frac{(\delta_1 x)^i}{i!} \ .$$

A recursive relationship exists among Laguerre polynomials:

$$\begin{aligned} L_0^{(a-1)}(-\delta_1 x) = 1, L_1^{(a-1)}(-\delta_1 x) = a+\delta_1 x,\ and \\ kL_k^{(a-1)}(-\delta_1 x) = (2k-2+a+\delta_1 x)L_{k-1}^{(a-1)}(-\delta_1 x)-(k+a-2)L_{k-2}^{(a-1)}(-\delta_1 x) \ . \end{aligned}$$

Laguerre polynomials of different orders can be obtained by recursive processes. However, the representation based on the Laguerre polynomial requires that $b$ is an integer, specifically that $n_2$ is an even number. This restricts the application scope of this method.

The abovementioned methods must calculate the sum of the infinite series. However, the actual calculation cannot proceed indefinitely. Thus, the truncating range must be determined for the numerical computation. The direct calculation method requires the pre-determination of the numbers of rows and columns involved in the calculation. And the index $j$ must be defined when using the method based on the Laguerre polynomial. If the preset truncating range is excessively small, the calculation result may be inaccurate. If the range is excessively large, a long computation time may be required. Therefore, in practical applications, a method that can determine the truncating range and conform to the precision requirement is required.

## 3. Proposed methods

The error between the real value and approximate value determines the calculating error of a numerical algorithm. However, the calculating error of Equation (4) is hard to obtain. Therefore, this paper will find an upper error bound that is slightly larger than the error, but easy to calculate. The upper error bound of doubly non-central beta or F distribution in numerical calculation has not been reported in the existing literature. However, many studies have focused on the upper error bound of the non-central F distribution. Norton [18] and Lenth [19] derived two upper bounds for non-central F distributions. Referring to the upper error bounds of Norton and Lenth and the division on the M-matrix, this paper establishes two methods named DIV1 and DIV2 to calculate the CDF of the doubly non-central beta distribution.

These two methods are enhanced variants of the direct calculation method, comparable in computational performance. And they can automatically determine the calculation range according to the precision requirement. Therefore, this paper recommends DIV1 and DIV2 as the numerical computation methods to calculate the CDF of the doubly non-central beta distribution.

### *3.1. DIV1*

DIV1 divides the **M**-matrix into three regions, as shown in Figure 2. Region 0 is the region involved in the numerical computation. Region 1 consists of the omitted items of each row in region 0. Region 2 consists of all the omitted rows below region 0. And we can assume that region 2 starts at row $j_1$. It is important to note that the boundaries between rows in Regions 1 and 0 may vary, meaning the boundary between these two regions does not correspond to a fixed column index.

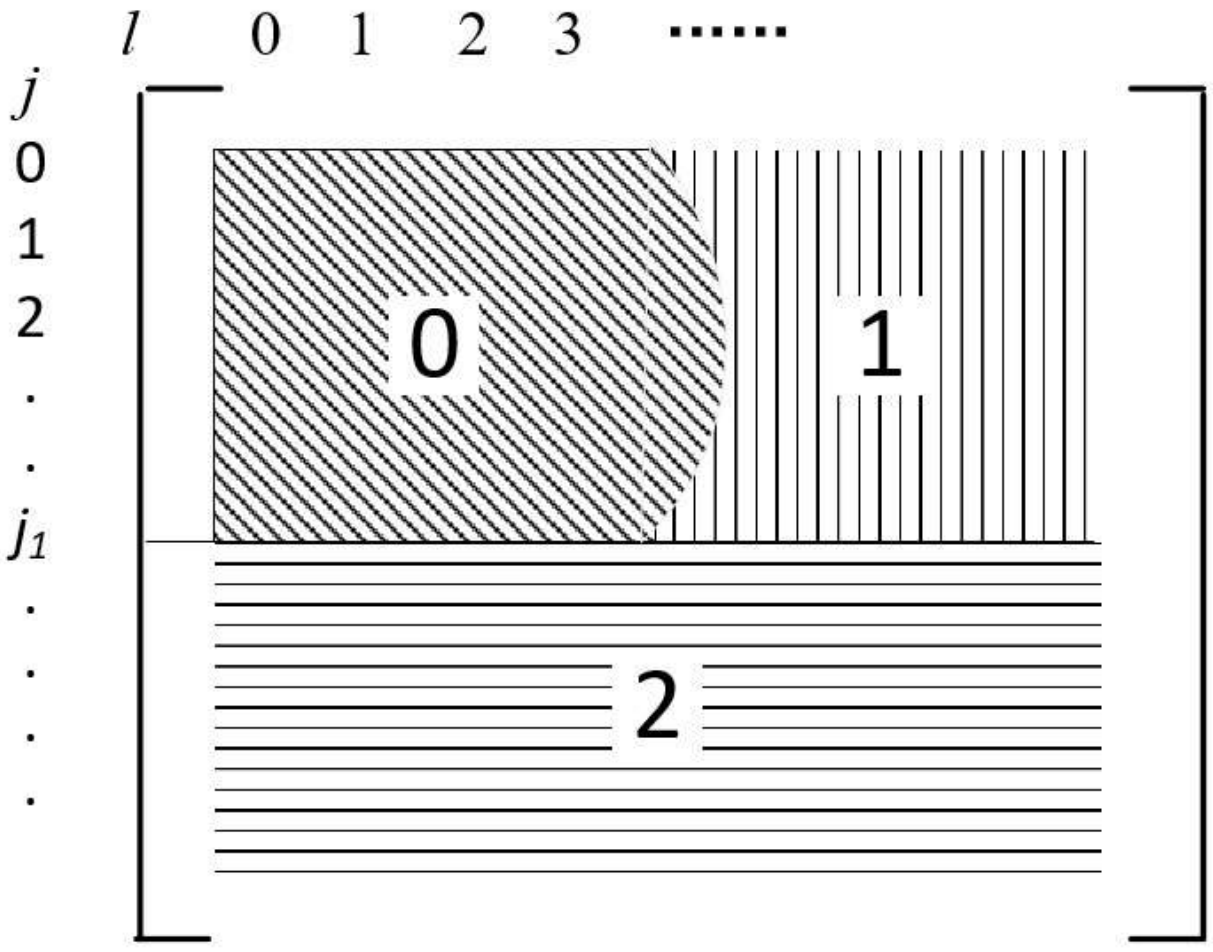

**Figure 2.** The division of the **M**-matrix in DIV1.

Let the sum of all items in region 0, 1 and 2 be $S_0$, $S_1$ and $S_2$, respectively. $S_0$ is the main component of calculating the CDF of the doubly non-central beta distribution. $S_1$+$S_2$ is the error of CDF of the doubly non-central beta distribution. $S_1$ and $S_2$ is difficult to compute, but the upper bounds of $S_1$ and $S_2$ can be determined.

Let $S_2$ be less than the presetting value $\varepsilon_{1,2}$, e.g., $10^{-5}$. Let $\varepsilon_{1,1}$ be the upper bound of the sum of all items in one row in region 1. If $j_1$ is determined, then $j_1\varepsilon_{1,1}$ is one upper bound of $S_1$, i.e., $S_1 < j_1\varepsilon_{1,1}$. It is easy to get $S_1$+$S_2 < j_1\varepsilon_{1,1} + \varepsilon_{1,\,2}$. The calculation precision of CDF can be controlled by $\varepsilon_{1,1}$ and $\varepsilon_{1,\,2}$. They are called the control lines in DIV1.

In this section, the upper bound of the sum of all items in region 2 is analyzed, and $j_1$ is obtained. Next, the upper bound of all items in every single row of region 1 is analyzed. Finally, the total upper error bound of DIV1 is established based on these two upper bounds.

#### 3.1.1. Upper bound in region 2

According to Theorem 2, and $e^x = \sum\limits_{n=0}^{\infty} \frac{x^n}{n!}$, the sum of all items in the $j^{th}$ row, $R_j(x : a, b; \delta_1, \delta_2)$, has the following relation:

$$\begin{aligned} R_j(x : a, b; \delta_1, \delta_2) &= e^{-(\delta_1+\delta_2)} \frac{\delta_1^j}{j!} \sum_{l=0}^{\infty} \frac{\delta_2^l}{l!} I_x(j+a, l+b) \\ &< e^{-(\delta_1+\delta_2)} \frac{\delta_1^j}{j!} e^{\delta_2} = e^{-\delta_1} \frac{\delta_1^j}{j!} \ . \end{aligned}$$

The sum of all the items following the $(j_1-1)^{th}$ row is

$$S_2 = R_{j_1} + R_{j_1+1} + \cdots < e^{-\delta_1}\left(\frac{\delta_1^{j_1}}{j_1!} + \frac{\delta_1^{j_1+1}}{(j_1+1)!} + \cdots\right)$$
$$= e^{-\delta_1}\left(e^{\delta_1} - \sum_{j=0}^{j_1-1}\frac{\delta_1^j}{j!}\right) = 1 - e^{-\delta_1}\sum_{j=0}^{j_1-1}\frac{\delta_1^j}{j!} \quad .$$

If $j_1$ is determined, then $1 - e^{-\delta_1}\sum_{j=0}^{j_1-1}\frac{\delta_1^j}{j!}$ is an upper bound in region 2. When $\varepsilon_{1,2}$ is known, $j_1$ can be obtained by solving the following inequality:

$$1 - e^{-\delta_1}\sum_{j=0}^{j_1-1}\frac{\delta_1^j}{j!} < \varepsilon_{1,2} \quad . \tag{8}$$

However, it is very difficult to find $j_1$ directly in inequality (8). Therefore, this paper searches $j$ from $0, 1$, to $\ldots$, until the first $j$ that makes the inequality (8) true .

#### 3.1.2. Upper bound in region 1

The series items in the **M**-matrix from row *0* to row $j_1-1$ , which are not included in region 0, form region 1. The sum in the $j^{th}$ row of the **M**-matrix truncated at $l = n_j - 1$ is as follows:

$$\hat{R}_j(x : a, b; \delta_1, \delta_2; n_j) = e^{-(\delta_1+\delta_2)}\frac{\delta_1^j}{j!}\sum_{l=0}^{n_j-1}\frac{\delta_2^l}{l!}I_x(j+a, l+b) \quad .$$

Let $e_j$ be the residual error of the $j^{th}$ row, defined as the difference between the theoretical value and the actual calculated value. In this case,

$$\begin{aligned} e_j =& R_j(x : a, b; \delta_1, \delta_2) - \hat{R}_j(x : a, b; \delta_1, \delta_2; n_j) \\ =& e^{-(\delta_1+\delta_2)}\frac{\delta_1^j}{j!}\sum_{l=n_j}^{\infty}\frac{\delta_2^l}{l!}I_x(j+a, l+b) \quad . \end{aligned} \tag{9}$$

Because $I_x(a, b) \le 1$, the following inequality can be obtained:

$$\sum_{l=n_j}^{\infty}\frac{\delta_2^l}{l!}I_x(j+a, l+b) \le \sum_{l=n_j}^{\infty}\frac{\delta_2^l}{l!} \quad .$$

Therefore, Equation (9) can be rewritten as

$$e_j < e^{-(\delta_1+\delta_2)}\frac{\delta_1^j}{j!}\sum_{l=n_j}^{\infty}\frac{\delta_2^l}{l!} = e^{-\delta_1}\frac{\delta_1^j}{j!}\left(1 - e^{-\delta_2}\sum_{l=0}^{n_j-1}\frac{\delta_2^l}{l!}\right) ,$$

where $0 \le j \le j_1 - 1$. If

$$e^{-\delta_1}\frac{\delta_1^j}{j!}\left(1 - e^{-\delta_2}\sum_{l=0}^{n_j-1}\frac{\delta_2^l}{l!}\right) < \varepsilon_{1,1} \ , \tag{10}$$

then the error in $j^{th}$ row can be controlled through $\varepsilon_{1,1}$. According to this derivation, the key to ensuring the calculation precision of each row in region 0 is to determine the truncated items $n_j$. However, $n_j$ cannot be directly obtained using inequality (10). Therefor, a numerical search must be performed for $l$ from 0 to $n_j - 1$. Obviously, in different rows, $n_j$ is different.

The upper bound in region 1 is as follows:

$$S_1 = e_0 + e_1 + \cdots + e_{j_1-1} < \sum_{j=0}^{j_1-1} e^{-\delta_1}\frac{\delta_1^j}{j!}\left(1 - e^{-\delta_2}\sum_{l=0}^{n_j-1}\frac{\delta_2^l}{l!}\right) < j_1\varepsilon_{1,1} \,. \tag{11}$$

3.1.3. Total upper error bound of DIV1

If the sum of all items after row $j_1 - 1$ in the **M**-matrix is less than $\varepsilon_{1,\,2}$, i.e., Equation (8) is true, and if the residual errors of row from 0 to $j_1 - 1$ are less than $\varepsilon_{1,1}$, i.e., Equation (10) is true; then the error between the true value $B_{a,\,b}^{\delta_1,\,\delta_2}(x)$ and the result of numerical calculation $\hat{B}_{a,\,b}^{\delta_1,\,\delta_2}(x)$ is

$$\begin{aligned} B_{a,\,b}^{\delta_1,\,\delta_2}(x) - \hat{B}_{a,\,b}^{\delta_1,\,\delta_2}(x) <& 1 - e^{-\delta_1}\sum_{j=0}^{j_1-1}\frac{\delta_1^j}{j!} + \sum_{j=0}^{j_1-1} e^{-\delta_1}\frac{\delta_1^j}{j!}\left(1 - e^{-\delta_2}\sum_{l=0}^{n_j-1}\frac{\delta_2^l}{l!}\right) \\ =& 1 - e^{-(\delta_1+\delta_2)}\sum_{j=0}^{j_1-1}\frac{\delta_1^j}{j!}\sum_{l=0}^{n_j-1}\frac{\delta_2^l}{l!} < j_1\varepsilon_{1,1} + \varepsilon_{1,\,2} \,. \end{aligned} \tag{12}$$

Let

$$U1_{a,b}^{\delta_1,\delta_2}(x) = 1 - e^{-(\delta_1+\delta_2)}\sum_{j=0}^{j_1-1}\frac{\delta_1^j}{j!}\sum_{l=0}^{n_j-1}\frac{\delta_2^l}{l!} \,, \tag{13}$$

then $U1_{a,b}^{\delta_1,\delta_2}(x)$ is the total upper error bound of DIV1. Obviously,

$$S_1 + S_2 < U1_{a,b}^{\delta_1,\delta_2}(x) < j_1\varepsilon_{1,1} + \varepsilon_{1,\,2} \,. \tag{14}$$

3.1.4. Pseudo-code of DIV1

The pseudo-codes of the DIV1 are shown in Algorithm 1. This algorithm involves two main steps. The first step is to determine the row boundary between regions 0 and 2, i.e., $j_1$, which are shown in line 2 to 7. The second step is to compute the sum of all items in one row that meet the precision requirement from row 0 to $j_1 - 1$, which are shown in line 9 to 20. The output is the probability value $P$.This algorithm considers the methods for avoiding numerical calculation errors as reported by [22].

### *3.2. DIV2*

Similar to DIV1, DIV2 divides the **M**-matrix into three regions, as shown in Figure 3. Region 0 is the region involved in the numerical computation. Region 1 consists of all the omitted columns after region 0, and we can assume that region 1 starts at column $l_1$. Region 2 consists of the omitted items of each column in region 0. The boundary between regions 2 and 0 is not a fixed row number. Let the sum of all items in region 0, 1 and 2 be $S_0$, $S_1$ and $S_2$, respectively. Let $S_1$ be less than the presetting value $\varepsilon_{2,1}$. Let $\varepsilon_{2,2}$ be the upper bound of the sum of all items in one column in region 2. $\varepsilon_{2,1}$ and $\varepsilon_{2,2}$ are called the control lines in DIV2. If $l_1$ is determined, then $l_1\varepsilon_{2,2}$ is one upper bound of $S_2$, i.e., $S_2 < l_1\varepsilon_{2,2}$. It is easy to get $S_1 + S_2 < \varepsilon_{2,1} + l_1\varepsilon_{2,2}$.

In this section, the upper bound of the sum of all items in region 1 is analyzed, and $l_1$ is obtained. Next, the upper bound of all items in every single column of region 2 is analyzed. Finally, the total upper error bound of DIV2 is determined on these two upper bounds.

3.2.1. Upper bound in region 1

Let $C_l(x : a, b; \delta_1, \delta_2)$ be the sum of all the items in the $l^{th}$ column. Because $I_x(a, b) \le 1$, the following inequality can be obtained:

**Algorithm 1** Pseudo-code of DIV1

**Input:** $x, a, b, \delta_1, \delta_2, \varepsilon_{1,1}$ and $\varepsilon_{1,2}$;
**Output:** $P$;
$t = 1,\ c = 1,\ j_1 = 0, S_2 = 1 - e^{-\delta_1}$;
**while** $S_2 > \varepsilon_{1,2}$ **do**
  $j_1 = j_1 + 1$;
  $t = t \cdot \delta_1 / j_1$;
  $c = c + t$;
  $S_2 = 1 - e^{-\delta_1} \cdot c$;
**end while**

$P = 0$;
**for** $j = 0,\ 1,\ 2, \cdots,\ j_1 - 1$ **do**
  $f = 0, l = 0, t1 = 1, tt = 1$;
  $ee = e^{-\delta_1} \delta_1^j / (j!)$;
  **while** $ee > \varepsilon_{1,1}$ **do**
    $f = f + e^{-(\delta_1 + \delta_2)} I_x(j + a, l + b) \delta_1^j \cdot t1 / (j!)$;
    $l = l + 1$;
    $t1 = t1 \cdot \delta_2 / l$;
    $tt = tt + t1$;
    $ee = e^{-\delta_1} \left(1 - e^{-\delta_2} \cdot tt\right) \delta_1^j / (j!)$;
  **end while**
  $P = P + e^{-(\delta_1 + \delta_2)} I_x(j + a, l + b) \delta_1^j \cdot t1 / (j!)$;
**end for**
**return** $P$.

$$
\begin{aligned}
C_l(x : a, b; \delta_1, \delta_2) &= e^{-(\delta_1 + \delta_2)} \frac{\delta_2^l}{l!} \sum_{j=0}^{\infty} \frac{\delta_1^j}{j!} I_x(j + a, l + b) \\
&< e^{-(\delta_1 + \delta_2)} \frac{\delta_2^l}{l!} \sum_{j=0}^{\infty} \frac{\delta_1^j}{j!} = e^{-\delta_2} \frac{\delta_2^l}{l!} \quad .
\end{aligned}
$$

Therefore, the sum of all the items after column $l_1 - 1$ is

$$
\begin{aligned}
S_1 &= C_{l_1} + C_{l_1+1} + \cdots \\
&< e^{-\delta_2} \left( \frac{\delta_2^{l_1}}{l_1!} + \frac{\delta_2^{l_1+1}}{(l_1+1)!} + \frac{\delta_2^{l_1+2}}{(l_1+2)!} \cdots \right) = 1 - e^{-\delta_2} \sum_{l=0}^{l_1-1} \frac{\delta_2^l}{l!} \quad .
\end{aligned}
$$

When $\varepsilon_{2,1}$ is known, $l_1$ can be obtained by solving the following inequality:

$$
1 - e^{-\delta_2} \sum_{l=0}^{l_1-1} \frac{\delta_2^l}{l!} < \varepsilon_{2,1} \quad . \tag{15}
$$

This paper determines $l_1$ by searcheing $l$ from 0, 1, to ... until the first $l$ that makes the inequality (15) true.

### 3.2.2. Upper bound in region 2

The sum in the $l^{th}$ column of **M**-matrix is truncated at row $j = m_l - 1$, as follows:

$$
\hat{C}_l(x : a, b; \delta_1, \delta_2; m_l) = e^{-(\delta_1 + \delta_2)} \frac{\delta_2^l}{l!} \sum_{j=0}^{m_l - 1} \frac{\delta_1^j}{j!} I_x(j + a, l + b) \quad .
$$

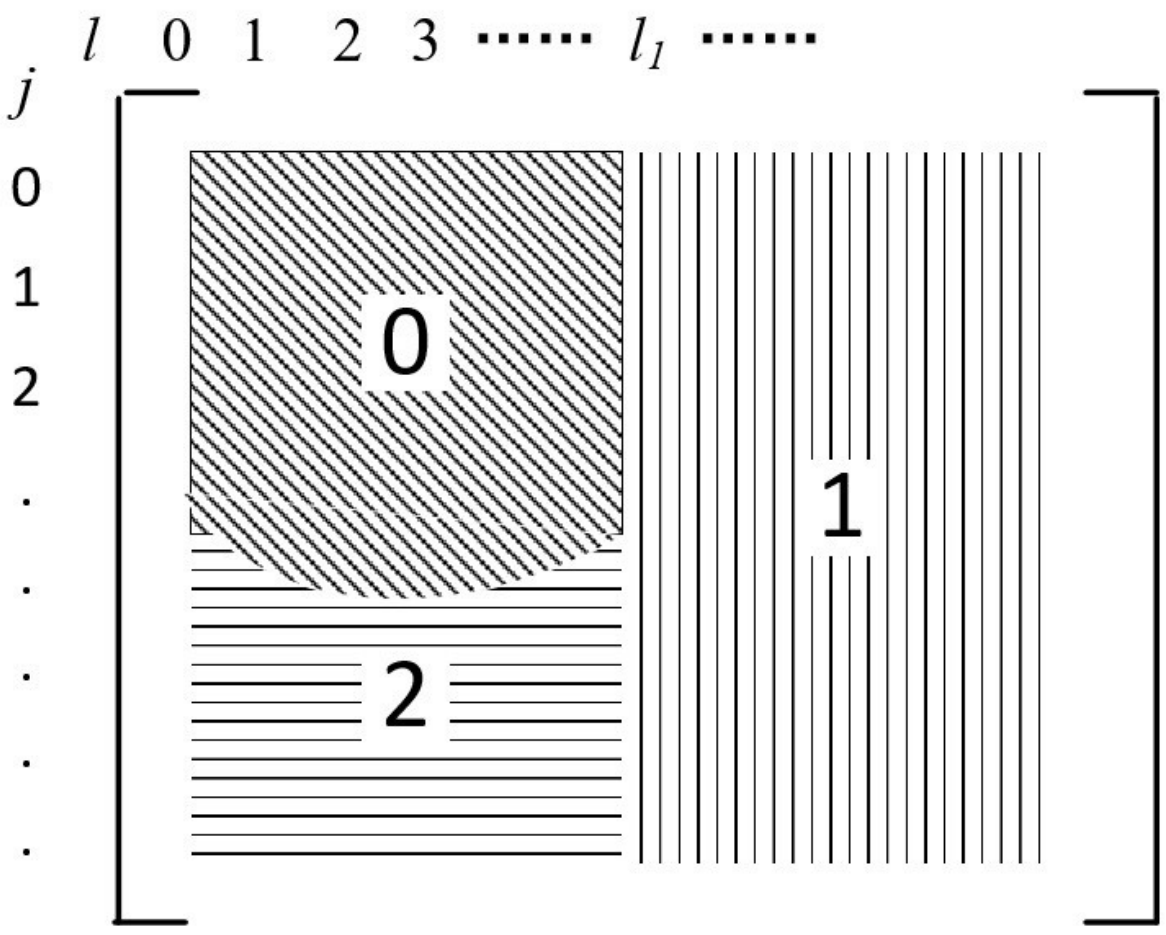

**Figure 3.** The division of the **M**-matrix in DIV2.

Let $e_l$ be the residual error of the $l^{th}$ row, defined as the difference between the theoretical value and the actual calculated value. In this case,

$$\begin{aligned} e_l &= C_l(x : a, b; \delta_1, \delta_2) - \hat{C}_l(x : a, b; \delta_1, \delta_2; m_l) \\ &= e^{-(\delta_1+\delta_2)} \frac{\delta_2^l}{l!} \sum_{j=m_l}^{\infty} \frac{\delta_1^j}{j!} I_x(j + a, l + b) \quad . \end{aligned} \tag{16}$$

Because $I_x(a, b) \leq 1$, then

$$\sum_{j=m_l}^{\infty} \frac{\delta_1^j}{j!} I_x(j + a, l + b) < \sum_{j=m_l}^{\infty} \frac{\delta_1^j}{j!} \quad . \tag{17}$$

Equation (16) can be rewritten as

$$\begin{aligned} e_l &< e^{-(\delta_1+\delta_2)} \frac{\delta_2^l}{l!} \sum_{j=m_l}^{\infty} \frac{\delta_1^j}{j!} \\ &= e^{-\delta_2} \frac{\delta_2^l}{l!} \left( 1 - e^{-\delta_1} \sum_{j=0}^{m_l-1} \frac{\delta_1^j}{j!} \right) \quad . \end{aligned}$$

If

$$e_l = e^{-\delta_2} \frac{\delta_2^l}{l!} \left( 1 - e^{-\delta_1} \sum_{j=0}^{m_l-1} \frac{\delta_1^j}{j!} \right) < \varepsilon_{2,2} \ , \tag{18}$$

the error in column $l$ can be controlled by $\varepsilon_{2,2}$. $m_l$ cannot be solved directly by using inequality (18), and a numerical search for $j$ must be performed from 0 to $m_l - 1$.

The upper bound in region 2 is as follows:

$$S_2 = e_0 + e_1 + \cdots + e_{l_1-1} < \sum_{l=0}^{l_1-1} e^{-\delta_2} \frac{\delta_2^l}{l!} \left( 1 - e^{-\delta_1} \sum_{j=0}^{m_l-1} \frac{\delta_1^j}{j!} \right) < l_1 \varepsilon_{2,2} \, . \tag{19}$$

### 3.2.3. Total upper error bound of DIV2

If the sumn of all items after column $l_1 - 1$ in the **M**-matrix is less than $\varepsilon_{2,1}$, i.e., Equation (15) is true, and if the residual errors of column from 0 to $l_1 - 1$ are less than $\varepsilon_{2,2}$, i.e., Equation (18) is true; then the error between the true value $B_{a,\,b}^{\delta_1,\,\delta_2}(x)$ and the result of numerical calculation $\hat{B}_{a,\,b}^{\delta_1,\,\delta_2}(x)$ is

$$\begin{aligned}
&B_{a,b}^{\delta_1,\delta_2}(x) - \hat{B}_{a,b}^{\delta_1,\delta_2}(x) \\
&< \left(1 - e^{-\delta_2}\sum_{l=0}^{l_1-1}\frac{\delta_2^l}{l!}\right) + \sum_{l=0}^{l_1-1} e^{-\delta_2}\frac{\delta_2^l}{l!}\left(1 - e^{-\delta_1}\sum_{j=0}^{m_l-1}\frac{\delta_1^j}{j!}\right) \\
&= 1 - e^{-(\delta_1+\delta_2)}\sum_{l=0}^{l_1-1}\frac{\delta_2^l}{l!}\sum_{j=0}^{m_l-1}\frac{\delta_1^j}{j!} < \varepsilon_{2,1} + l_1\varepsilon_{2,2}\ .
\end{aligned}$$

Let

$$U2_{a,b}^{\delta_1,\delta_2}(x) = 1 - e^{-(\delta_1+\delta_2)}\sum_{l=0}^{l_1-1}\frac{\delta_2^l}{l!}\sum_{j=0}^{m_l-1}\frac{\delta_1^j}{j!}\ , \tag{20}$$

then $U2_{a,b}^{\delta_1,\delta_2}(x)$ is the total error upper bound of DIV2. Obviously,

$$S_1 + S_2 < U2_{a,b}^{\delta_1,\delta_2}(x) < \varepsilon_{2,1} + l_1\varepsilon_{2,2}\ . \tag{21}$$

#### 3.2.4. Pseudo-code of DIV2

The pseudo-codes of the DIV2 are shown in Algorithm 2. This algorithm involves two main steps. The first step is to determine the column boundary $l_1$ between regions 0 and 1, which are shown in line 2 to 7. The second step is to compute the sum of all items in one column that satisfy the precision requirement from column 0 to $l_1 - 1$, which are shown in line 9 to 20. The output is the probability value $P$.

**Algorithm 2** Pseudo-code of DIV2

**Input:** $x, a, b, \delta_1, \delta_2, \varepsilon_{2,1}$ and $\varepsilon_{2,2}$;
**Output:** $P$;

1: $t = 1,\ c = 1,\ l_1 = 0, S_1 = 1 - e^{-\delta_2}$;
2: **while** $S_1 > \varepsilon_{2,1}$ **do**
3: $\quad l_1 = l_1 + 1$;
4: $\quad t = t \cdot \delta_2 / l_1$;
5: $\quad c = c + t$;
6: $\quad S_2 = 1 - e^{-\delta_2} \cdot c$;
7: **end while**
8:
9: $P = 0$;
10: **for** $l = 0,\ 1,\ 2, \cdots,\ l_1 - 1$ **do**
11: $\quad f = 0, j = 0, t1 = 1, tt = 1$;
12: $\quad ee = e^{-\delta_2}\delta_2^l/(l!)$;
13: $\quad$**while** $ee > \varepsilon_{2,2}$ **do**
14: $\qquad f = f + e^{-(\delta_1+\delta_2)} I_x(j+a, l+b) \cdot t1\delta_2^l/(l!)$;
15: $\qquad j = j + 1$;
16: $\qquad t1 = t1 \cdot \delta_1 / j$;
17: $\qquad tt = tt + t1$;
18: $\qquad ee = e^{-\delta_2}(1 - e^{-\delta_1} \cdot tt)\delta_2^l/(l!)$;
19: $\quad$**end while**
20: $\quad P = P + e^{-(\delta_1+\delta_2)} I_x(j+a, l+b) \cdot t1\delta_2^l/(l!)$;
21: **end for**
22: **return** $P$.

## 4. Numerical experiments

Considering different parameter states, the computation precision and speed of the direct calculation method were compared with three proposed methods.

In the experiments, if no hint is given, then $\varepsilon_{1,1} = \varepsilon_{2,2} = 10^{-5}$ and $\varepsilon_{1,2} = \varepsilon_{2,1} = 10^{-7}$. To obtain accurate calculation results and facilitate comparison with the proposed methods, the direct calculation method was used to accurately calculate the sum of the rows and columns and the entire M-matrix. In the direct calculation, we truncated at least 100 items in each row or column to ensure that the calculation result was accurate to 8 decimal places.

All the experiments were conducted on a laptop with a Core i5 CPU, 16 GB memory and 200 GB SSD. All the programmes were implemented in MATLAB 2021a.

### *4.1. Precision*

#### 4.1.1. Total error

Under different parameters, the comparison of the results and precision obtained using the direct calculation method and three methods presented in this paper are summarized in Table 1. *P0* is the probability determined in the direct calculation. *P1* and *P2* are the probabilities calculated using DIV1 and DIV2, respectively. *Error1*,and *Error2* indicate the differences between the exact calculation results and the results obtained using two proposed methods. *UB1* and *UB2* are the total upper error bounds presented in Equation (13) and (20), respectively. *CL1* and *CL2*, determined by $\varepsilon_{1,1}$, $\varepsilon_{1,2}$, $\varepsilon_{2,1}$ and $\varepsilon_{2,2}$, represent the total control lines of the error respectively for two division methods. $j_1$ and $l_1$ are the row boundaries of region 1 in DIV1 and column boundaries in DIV2, respectively.

It can be noted that the calculation errors *Error1* are always less than the total upper error bounds *UB1*, and the total upper error bounds*UB1* are always less than the total control lines *CL1*. The same relationship holds for *Error2*, *UB2* and *CL2*. In other words, the total upper error bounds determined by the total control lines in this paper is higher than the actual calculation error, which indicates that the upper error bound obtained using the proposed methods is reasonable. The upper error bound can effectively control the calculation precision of the CDF of the doubly non-central beta distribution.

Two division methods can also calculate the probability of a doubly non-central F distribution through Equation (7). The comparison results with Tiku's and DIV1 method are shown in Table 2. *P0*, *DIV1* and *Tiku* are the probability determined in the direct calculation, DIV1 and Tiku's methods, respectively. *Error1* and *Error2* indicate the differences between the exact probabilities and the probabilities obtained using DIV1 and Tiku's methods. The probabilities in penultimate column of Table 2 is the results selected from the Table 3 and 4 in [16], which is only up to four decimal places.

It can be seen that the absolute values of *Error1* are all less than those of *Error2*. Usually, the probabilities calculated by Tiku's method are accurate to three decimal places, but sometimes it is only two decimal places, as seen in line 4 of Table 2. DIV1 can calculate the result to meet the preset precision, while Tiku's method cannot.

#### 4.1.2. Column and row error

In this experiment, the column boundaries between regions 0 and 1 and calculation errors of each row in region 0 for DIV1 were obtained, as shown in Table 3. Moreover, the row boundaries between regions 0 and 2 and calculation errors of each column in region 0 for DIV2 were obtained, as shown in Table 4.

The parameters in this experiment were set as $n_1 = 5$, $n_2 = 7$, $\delta_1 = 6.25$, $\delta_2 = 0.25$ and $x = 0.3$. In Table 3, $j$ is the row index, $n_j$ is the number of truncated items in a row, $\hat{R}_j$ is the sum of all the truncated items in a row, $R_j$ is the result of the exact calculation, and $e_j$ is the value of $R_j$ minus $\hat{R}_j$, representing the calculation error. $UB_e_j$ calculated by Equation (10) is the upper error bound in the $j^{th}$ row.

In Table 4, $l$ is the column index, $m_l$ is the number of truncated items in a column, $\hat{C}_l$ is the sum of all the truncated items in a column, $C_l$ is the result of the exact calculation, and $e_l$ is the value of $C_l$ minus $\hat{C}_l$. $UB_e_l$ calculated by Equation (18) is the upper error bound in the $l^{th}$ column.

**Table 1.** Comparison of precision associated with the direct calculation method and two proposed methods under different parameters.

| | | | | | | DIV1 | | | | DIV2 | | | |
|---|---|---|---|---|---|---|---|---|---|---|---|---|---|
| $n_1$ | $n_2$ | $\lambda_1$ | $\lambda_2$ | $x$ | P0(Exact) | P1 | Error1 | UB1 | CL1 | P2 | Error2 | UB2 | CL2 |
| 2 | 4 | 0.5 | 0.5 | 0.7 | 0.8967439 | 0.8967413 | $2.61\times10^{-6}$ | $6.71\times10^{-6}$ | $1.05\times10^{-5}$ | 0.8967373 | $6.67\times10^{-6}$ | $6.71\times10^{-6}$ | $1.05\times10^{-5}$ |
| 3 | 6 | 1 | 2 | 0.3 | 0.4843354 | 0.4843352 | $2.03\times10^{-7}$ | $1.36\times10^{-6}$ | $1.07\times10^{-5}$ | 0.4843343 | $1.07\times10^{-6}$ | $1.42\times10^{-6}$ | $1.09\times10^{-5}$ |
| 4 | 30 | 24 | 5 | 0.8 | 0.9999335 | 0.9999232 | $1.02\times10^{-5}$ | $1.03\times10^{-5}$ | $1.30\times10^{-5}$ | 0.9999302 | $3.28\times10^{-6}$ | $3.29\times10^{-6}$ | $1.13\times10^{-5}$ |
| 30 | 4 | 24 | 5 | 0.8 | 0.2114543 | 0.2114528 | $1.53\times10^{-6}$ | $1.03\times10^{-5}$ | $1.30\times10^{-5}$ | 0.2114517 | $2.60\times10^{-6}$ | $3.29\times10^{-6}$ | $1.13\times10^{-5}$ |
| 5 | 7 | 0.25 | 6.25 | 0.3 | 0.5685838 | 0.5685829 | $8.85\times10^{-7}$ | $9.33\times10^{-6}$ | $1.04\times10^{-5}$ | 0.5685784 | $5.35\times10^{-6}$ | $5.85\times10^{-6}$ | $1.14\times10^{-5}$ |
| 5 | 7 | 6.25 | 0.25 | 0.3 | 0.0593471 | 0.0593471 | $6.64\times10^{-8}$ | $5.85\times10^{-6}$ | $1.14\times10^{-5}$ | 0.0593450 | $2.14\times10^{-6}$ | $9.33\times10^{-6}$ | $1.04\times10^{-5}$ |
| 6 | 8 | 5 | 25 | 0.3 | 0.6877595 | 0.6877587 | $8.42\times10^{-7}$ | $3.37\times10^{-6}$ | $1.13\times10^{-5}$ | 0.6877519 | $7.60\times10^{-6}$ | $8.98\times10^{-6}$ | $1.31\times10^{-5}$ |
| 8 | 15 | 4 | 9 | 0.6 | 0.9756436 | 0.9756376 | $5.97\times10^{-6}$ | $8.98\times10^{-6}$ | $1.11\times10^{-5}$ | 0.9756377 | $5.87\times10^{-6}$ | $6.09\times10^{-6}$ | $1.17\times10^{-5}$ |

**Table 2.** Precision comparison for CDF of doubly non-central F distribution with Tiku and DIV1 method under different parameters.

| $n_1$ | $n_2$ | $\lambda_1$ | $\lambda_2$ | $f$ | $P0$ | $DIV1$ | $Error1$ | $Tiku$ | $Error2$ |
|---|---|---|---|---|---|---|---|---|---|
| 2 | 4 | 1.5 | 1.5 | 6.94414 | 0.933730 | 0.933729 | $9.38\times10^{-7}$ | 0.9325 | $1.23\times10^{-3}$ |
| 2 | 15 | 1.5 | 3 | 3.68235 | 0.893163 | 0.893163 | $3.14\times10^{-7}$ | 0.8898 | $3.36\times10^{-3}$ |
| 4 | 30 | 2 | 2 | 2.68966 | 0.871013 | 0.871013 | $3.07\times10^{-7}$ | 0.8704 | $6.13\times10^{-4}$ |
| 8 | 15 | 4 | 9 | 2.64079 | 0.968629 | 0.968623 | $5.52\times10^{-6}$ | 0.9415 | $2.71\times10^{-2}$ |
| 2 | 4 | 12 | 3 | 6.94414 | 0.711489 | 0.711487 | $1.80\times10^{-6}$ | 0.7138 | $-2.31\times10^{-3}$ |
| 4 | 30 | 24 | 5 | 2.68966 | 0.057048 | 0.057047 | $4.48\times10^{-7}$ | 0.0513 | $5.75\times10^{-3}$ |

**Table 3.** Row calculation results of region 0 in DIV1.

| $j$ | $n_j$ | $R_j$ | $\hat{R}_j$ | $e_j$ | $UB_e_j$ |
|---|---|---|---|---|---|
| 0 | 4 | 0.013639 | 0.013639 | $7.51\times10^{-9}$ | $1.01\times10^{-8}$ |
| 1 | 4 | 0.021005 | 0.021005 | $1.77\times10^{-8}$ | $3.15\times10^{-8}$ |
| 2 | 4 | 0.015023 | 0.015023 | $1.90\times10^{-8}$ | $4.92\times10^{-8}$ |
| 3 | 4 | 0.006785 | 0.006785 | $1.27\times10^{-8}$ | $5.12\times10^{-8}$ |
| 4 | 4 | 0.002205 | 0.002205 | $5.93\times10^{-9}$ | $4.00\times10^{-8}$ |
| 5 | 4 | 0.000555 | 0.000555 | $2.10\times10^{-9}$ | $2.50\times10^{-8}$ |
| 6 | 4 | 0.000113 | 0.000113 | $5.90\times10^{-10}$ | $1.30\times10^{-8}$ |
| 7 | 4 | $1.94\times10^{-5}$ | $1.94\times10^{-5}$ | $1.37\times10^{-10}$ | $5.81\times10^{-9}$ |
| 8 | 3 | $2.86\times10^{-6}$ | $2.86\times10^{-6}$ | $6.13\times10^{-10}$ | $9.12\times10^{-8}$ |
| 9 | 3 | $3.68\times10^{-7}$ | $3.68\times10^{-7}$ | $9.80\times10^{-11}$ | $3.17\times10^{-8}$ |
| 10 | 3 | $4.21\times10^{-8}$ | $4.21\times10^{-8}$ | $1.37\times10^{-11}$ | $9.90\times10^{-9}$ |
| 11 | 2 | $4.33\times10^{-9}$ | $4.30\times10^{-9}$ | $2.51\times10^{-11}$ | $9.06\times10^{-8}$ |
| 12 | 2 | $4.03\times10^{-10}$ | $4.01\times10^{-10}$ | $2.68\times10^{-12}$ | $2.36\times10^{-8}$ |
| 13 | 2 | $3.44\times10^{-11}$ | $3.42\times10^{-11}$ | $2.60\times10^{-13}$ | $5.67\times10^{-9}$ |

In Table 3, the calculation error $e_j$ is always less than the upper error bounds in $j^th$ row $UB_e_j$, and the upper error bound in $j^th$ row $UB_e_j$ is always less than the row control line $\varepsilon_{1,1}$. The same relationship holds for $e_l$, $UB_e_l$ and $\varepsilon_{2,2}$. This finding indicates that the proposed methods can control the calculation precision of each row or column according to the error bound.

In the second column of Table 3, the numbers of the truncated items in different rows are different, which indicates that the column boundaries of regions 0 and 1 are not a fixed column index in DIV1. Similarly, the row boundaries of regions 0 and 2 in DIV2 are not a fixed row index.

#### 4.1.3. Precision control

In DIV1 and DIV2, the role of parameters $\varepsilon_{1,1}$, $\varepsilon_{1,2}$, $\varepsilon_{2,1}$, and $\varepsilon_{2,2}$ is to control the precision of the algorithms. The effects of different parameter settings on precision are shown in Table 5, where $n_1$, $n_2$, $\lambda_1$, $\lambda_2$, $x$ are 8, 15, 4, 9, 0.6, respectively. Since the precision of DIV1 and DIV2 is similar, only DIV1 is chosen in this experiment. *Error1* is the differences between the exact calculation results and DIV1. *UB1* is the total upper error bounds. *CL1* is the total control lines of the error.

**Table 4.** Column calculation results of region 0 in DIV2.

| $l$ | $m_l$ | $C_l$ | $\hat{C}_l$ | $e_l$ | $UB_e_l$ |
|---|---|---|---|---|---|
| 0 | 16 | 0.048552 | 0.048552 | $5.27\times10^{-16}$ | $3.40\times10^{-8}$ |
| 1 | 15 | 0.009868 | 0.009868 | $4.70\times10^{-15}$ | $2.34\times10^{-8}$ |
| 2 | 13 | 0.000904 | 0.000904 | $1.93\times10^{-13}$ | $3.70\times10^{-8}$ |
| 3 | 11 | $5.13\times10^{-5}$ | $5.13\times10^{-5}$ | $2.75\times10^{-12}$ | $2.99\times10^{-8}$ |

In Table 5, when $\varepsilon_{1,1}$ and $\varepsilon_{1,2}$ are smaller, the calculation precision is higher. $\varepsilon_{1,2}$ has a greater effect on precision than $\varepsilon_{1,1}$. This experiment shows that we can make the algorithm achieve any precision by controlling $\varepsilon_{1,2}$ and $\varepsilon_{1,1}$.

**Table 5.** Influence of control parameters on calculation precision.

| $\varepsilon_{1,2}$ | $\varepsilon_{1,1}$ | $P$ | *Error1* | *UB1* | *CL1* |
|---|---|---|---|---|---|
| $1.0\times10^{-3}$ | $1.0\times10^{-5}$ | 0.975403 | $2.41\times10^{-4}$ | $3.01\times10^{-4}$ | $1.09\times10^{-3}$ |
| $1.0\times10^{-4}$ | $1.0\times10^{-5}$ | 0.975540 | $1.04\times10^{-4}$ | $1.18\times10^{-4}$ | $2.00\times10^{-3}$ |
| $1.0\times10^{-4}$ | $1.0\times10^{-6}$ | 0.975605 | $3.87\times10^{-5}$ | $5.30\times10^{-5}$ | $1.10\times10^{-4}$ |
| $1.0\times10^{-5}$ | $1.0\times10^{-6}$ | 0.975631 | $1.24\times10^{-5}$ | $1.54\times10^{-5}$ | $2.10\times10^{-5}$ |
| $1.0\times10^{-5}$ | $1.0\times10^{-7}$ | 0.975638 | $5.97\times10^{-6}$ | $8.98\times10^{-6}$ | $1.11\times10^{-5}$ |
| $1.0\times10^{-6}$ | $1.0\times10^{-7}$ | 0.975643 | $8.65\times10^{-7}$ | $9.70\times10^{-7}$ | $2.30\times10^{-6}$ |
| $1.0\times10^{-6}$ | $1.0\times10^{-8}$ | 0.975643 | $1.76\times10^{-7}$ | $2.75\times10^{-7}$ | $1.13\times10^{-6}$ |

### 4.2. Computational speed

The comparison of the computational speed of the direct calculation and the two proposed methods is shown in Table 6. In this table, *Exact* refers to the direct calculation, *ItemNum* is the count of items involved in the calculation of the two division methods. *Time(sec)* is the computational time, and the unit is second. This paper truncated 100 items in each row or column to ensure that the direct calculation result was precise. The experiment was conducted 80 times, and the average value was used for *Time(sec)*.

**Table 6.** Comparison of computational speed of the direct calculation and two proposed methods.

| | | | | | *Exact* | *DIV1* | | *DIV2* | |
|---|---|---|---|---|---|---|---|---|---|
| $n_1$ | $n_2$ | $\lambda_1$ | $\lambda_2$ | $x$ | *Time(sec)* | *itemNum* | *Time(sec)* | *itemNum* | *Time(sec)* |
| 2 | 4 | 0.5 | 0.5 | 0.7 | $6.89\times10^{-3}$ | 28 | $1.00\times10^{-4}$ | 28 | $1.03\times10^{-4}$ |
| 3 | 6 | 1 | 2 | 0.3 | $7.29\times10^{-3}$ | 57 | $7.04\times10^{-5}$ | 58 | $7.67\times10^{-5}$ |
| 4 | 30 | 24 | 5 | 0.8 | $6.53\times10^{-3}$ | 362 | $2.53\times10^{-4}$ | 388 | $2.68\times10^{-4}$ |
| 30 | 4 | 24 | 5 | 0.8 | $6.77\times10^{-3}$ | 362 | $2.65\times10^{-4}$ | 388 | $2.64\times10^{-4}$ |
| 5 | 7 | 0.25 | 6.25 | 0.3 | $7.25\times10^{-3}$ | 59 | $7.16\times10^{-5}$ | 61 | $6.79\times10^{-5}$ |
| 5 | 7 | 6.25 | 0.25 | 0.3 | $7.21\times10^{-3}$ | 61 | $5.20\times10^{-5}$ | 59 | $4.42\times10^{-5}$ |
| 6 | 8 | 5 | 25 | 0.3 | $6.91\times10^{-3}$ | 397 | $2.65\times10^{-4}$ | 371 | $2.57\times10^{-4}$ |
| 8 | 15 | 4 | 9 | 0.6 | $7.47\times10^{-3}$ | 187 | $1.41\times10^{-4}$ | 191 | $1.41\times10^{-4}$ |

According to the experimental results, the computational times of DIV1 and DIV2 are comparable, and both are less than that of the direct calculation method. Because the number of items involved in DIV1 and DIV2 are less than those in the direct calculation.

## 5. Discussions

(1) Preset values for $\varepsilon_{1,1}$, $\varepsilon_{1,2}$, $\varepsilon_{2,1}$, and $\varepsilon_{2,2}$. In the computation of DIV1, the number of items truncated in one row rarely exceeded 100. Therefore, $\varepsilon_{1,1}$ is usually set as 1% of $\varepsilon_{1,2}$, which can ensure that $\varepsilon_{1,1}$ in the upper bound of the total error is less than $\varepsilon_{1,2}$. For the same reason, $\varepsilon_{2,2}$ is usually set as 1% of $\varepsilon_{2,1}$.

(2) Control lines can be used instead of upper error bounds. Two upper error bounds, DIV1 and DIV2, are difficult to calculate. but two control Lines, $j_1\varepsilon_{1,1} + \varepsilon_{1,2}$ and $\varepsilon_{2,1} + l_1\varepsilon_{2,2}$, are easy to calculate. Therefore, we recommend using control lines instead of upper error bounds in practice.

(3) Accelerating the computation of the incomplete beta function. Baharev [22] showed that incomplete beta functions can be computed in a recursive manner. However, this paper used the betainc function provided by MATLAB and did not use the recursive method. To enhance the algorithm, the recursive method can be used to increase the computational speed.

(4) According to the experimental results, the upper error bound determined in this paper is considerably higher than the actual error in certain cases. As shown in Line 5 of Table 1, the upper error bound of the two methods is 10 times the actual error. Future work can focus on examining whether an upper bound similar to the actual error value can be obtained.

(5) Although the two methods presented in this paper can ensure the calculation precision, their computational speeds are slower than Tiku's method.

## 6. Conclusion

In this paper, we propose two numerical computing methods for the CDF of the doubly non-central beta distribution based on the M-matrix division. We derive the theoretical upper error bound of the two methods and present the computation steps through pseudo-codes. Both methods can automatically compute the values after setting the control lines. It is not necessary to set the calculation range of the M-matrix, and the calculation precision can be calculated. Therefore, we recommend the two methods as the numerical computation methods for the CDF of the doubly non-central beta distribution.

**Author Contributions:** Conceptualization, J.W. and H.L.; methodology, H.L.; software, H.L. and L.F.; validation, J.W. and F.M.; formal analysis, F.M.; writing—original draft preparation, Y.T.; writing—review and editing, T.D.; funding acquisition, Y.T. All authors have read and agreed to the published version of the manuscript.

**Funding:** This research was funded by the research project funded by the National Natural Science Foundation of China grant number 72101137, the Natural Science Foundation of Shandong Province grant number ZR2021MF117, the Education Ministry Humanities and Social Science Research Youth Fund Project of China grant number 21YJCZH150.

**Conflicts of Interest:** The authors declare no conflicts of interest.